\documentclass[a4paper, 11pt]{article}

\pdfoutput=1 % makes arXiv use pdflatex

\usepackage[utf8]{inputenc}
\usepackage[textwidth=5in,textheight=8.5in]{geometry}
\usepackage[english]{babel}
\usepackage{amsmath}
\usepackage{amssymb}
\usepackage{amsthm} 
\usepackage[alwaysadjust]{paralist}
\usepackage[%
breaklinks,
colorlinks=true,
linkcolor=black,
anchorcolor=black,
citecolor=black,
filecolor=black,
menucolor=black,
urlcolor=black]{hyperref}
\usepackage{microtype}

\newcommand{\R}{\mathbb{R}}
\newcommand{\Z}{\mathbb{Z}}
\newcommand{\RP}{\R\mathrm{P}}
\newcommand{\Moeb}{\textrm{Möb}}
\newcommand{\Sim}{\operatorname{Sim}}
\newcommand{\PO}{\operatorname{PO}}
\newcommand{\st}{\operatorname{star}}

\theoremstyle{definition}
\newtheorem{definition}{Definition}[section]

\newtheorem{remark}[definition]{Remark}

\theoremstyle{plain}
\newtheorem{theorem}{Theorem}
\newtheorem*{theorem*}{Theorem}

\newtheorem{lemma}[definition]{Lemma}

\title{A discrete version of Liouville's theorem on conformal maps}

\author{Ulrich Pinkall, Boris Springborn}

\date{}

\begin{document}

\maketitle

\begin{abstract}
  Liouville's theorem says that in dimension greater than two, all
  conformal maps are Möbius transformations. We prove an analogous
  statement about simplicial complexes, where two simplicial complexes
  are considered discretely conformally equivalent if they are
  combinatorially equivalent and the lengths of corresponding edges
  are related by scale factors associated with the vertices.
  
  \medskip\noindent%

  \medskip\noindent%
\end{abstract}

\section{Introduction}
\label{sec:intro}

Liouville's theorem says that in dimension three and higher, conformal
maps are Möbius transformations. More precisely:

\begin{theorem*}[Liouville]
  If $U\subset\R^{n}$ is a domain and $n\geq 3$, then any sufficiently
  regular conformal map $f:U\rightarrow\R^{n}$ is the restriction of a
  Möbius transformation.
\end{theorem*}

Liouville himself proved this theorem under the assumption that $f$ is
four times continuously differentiable~\cite[Note VI,
pp.~609--616]{monge_liouville_1850}. Different modern proofs can be
found, e.g., in the textbooks of Blaschke~\cite[\S40,
pp.~66f]{blaschke_1921} and Dubrovin, Fomenko \& Novikov~\cite[\S15,
pp.~138ff]{dfn_1992}. The regularity assumption can be weakened
considerably. Gehring~\cite{gehring_1962} and
Reshetnyak~\cite{reshetnyak_1967} showed it is sufficient to assume
that $f$ is in the Sobolev space $W^{1}_{n,\textit{loc}}(U)$, see
also~\cite{bojarski82} and~\cite{iwaniec_martin_1993}.

The purpose of this article is to extend the theorem in a different
direction. We establish the following discrete version of Liouville's
theorem for simplicial complexes (see Section~\ref{sec:basic} for
precise definitions):

\begin{theorem}
  \label{thm:dliouville}
  If $n\geq 3$, then two locally Delaunay discrete domains in $\R^{n}$
  are discretely conformally equivalent if and only if they are
  Möbius equivalent.
\end{theorem}

Roughly, simplicial complexes are considered \emph{discretely
  conformally equivalent} if they are combinatorially equivalent and
the lengths of corresponding edges are related by scale factors
associated with the vertices. They are considered \emph{Möbius
  equivalent} if they are combinatorially equivalent and the vertex
positions are related by a Möbius transformation.

One implication of the equivalence statement, ``Möbius equivalence
implies discrete conformal equivalence'', holds for arbitrary
simplicial complexes and for any dimension~$n$ (see
Section~\ref{sec:easy}). The other implication, ``discrete conformal
equivalence implies Möbius equivalence'', is only true for $n\geq 3$
and for a more restrictive class of simplicial complexes (see Section~\ref{sec:harder}). In
Definitions~\ref{def:discrete_domain} and~\ref{def:Delaunay} we
therefore define a \emph{locally Delaunay discrete domain} to be a
locally finite, full-dimensional simplicial complex that satisfies
some additional conditions that are sufficient and necessary to deduce
Möbius equivalence from discrete conformal equivalence in dimension
three or greater.

The basic concepts are explained in
Section~\ref{sec:basic}. Section~\ref{sec:proof} is devoted to a proof
of Theorem~\ref{thm:dliouville}, which may in hindsight appear rather
obvious. Theorem~\ref{thm:dliouville} and its proof suggest a
necessary and sufficient condition for discrete conformal
flatness. This, the connection between discrete conformal equivalence
and hyperbolic geometry, and some open questions are discussed
in Section~\ref{sec:outlook}.

A related approach to discretize the notion of conformality is via
circle packings~\cite{Stephenson}, or, in higher dimension, sphere
packings. Cooper~\& Rivin's local rigidity theorem for sphere
packings~\cite{Cooper-Rivin, Glickenstein, Rivin} and Xu's global
version~\cite{Xu} state, roughly, that a 3-manifold triangulation and
prescribed solid angles at the vertices determine a sphere packing
uniquely (if it exists). An analogous statement may also hold in the
setting of discrete conformal equivalence, roughly: A 3-manifold
triangulation together with functions that assign a length to each
edge and an solid angle to each vertex determine a discretely
conformally equivalent triangulation with the given solid angles at
the vertices uniquely (if it exists). It seems natural to expect that
the analytic method of proof, based on a variational principle,
also extends to this setting.

But our Theorem~\ref{thm:dliouville} is different in nature: It is not
about the metric rigidity of piecewise flat closed manifolds, but
about the Möbius rigidity of triangulated domains in $\R^{n}$. The
method of proof is also very different: Rather than a variational
principle, our proof of Theorem~\ref{thm:dliouville} relies on
Cauchy's rigidity theorem for convex polyhedra and its higher
dimensional generalization~\cite{Pak}, which does all the hard work
(see Section~\ref{sec:harder}). The rest is essentially just setup.
Cauchy's rigidity theorem will be applied to Möbius images of the
links of interior vertices, with convexity ensured by the local
Delaunay condition. This explains why the discrete Liouville theorem
holds only in dimensions three and higher: In dimension two, the
respective links are polygons, which are not rigid.
 
\section{Basic definitions}
\label{sec:basic}

In this article, a \emph{combinatorial isomorphism} of simplicial
complexes~$K$ and~$K'$ in~$\R^{n}$ is understood to be a bijection
\begin{equation}
  \label{eq:phi}
  \phi:V\longrightarrow V'
\end{equation}
between the vertex sets $V$ and $V'$ of $K$ and $K'$, respectively,
such that for any subset $\{v_{0},\ldots,v_{k}\}\subseteq V$ the
simplex
\begin{equation*}
  [v_{0},\ldots,v_{k}]\subset\R^{n}
\end{equation*}
is an element of $K$ if and only if the simplex
\begin{equation*}
  [\phi(v_{1}),\ldots,\phi(v_{k})]\subset\R^{n}
\end{equation*}
is an element of $K'$. Thus, a combinatorial isomorphism $\phi$ induces a
bijection between the complexes $K$ and $K'$, as well as a piecewise
linear simplicial map between their carriers $|K|$ and
$|K'|$. Simplicial complexes are \emph{combinatorially equivalent} if
there exists a combinatorial isomorphism between their vertex sets.

\begin{definition}[discrete conformal equivalence]
  \label{def:dce}
  Combinatorially equivalent simplicial complexes $K$ and $K'$ in
  $\R^{n}$ are called \emph{discretely conformally equivalent with
    respect to a combinatorial isomorphism}~\eqref{eq:phi} if there
  exists a function
  \begin{equation*}
    u:V\longrightarrow\R
  \end{equation*}
  such that the length of each edge $[v_{1},v_{2}]\in K$ is related to
  the length of the corresponding edge
  $[\phi(v_{1}),\phi(v_{2})]\in K'$, by
  \begin{equation}
    \label{eq:scale}
    \big|\,\phi(v_{1})\,-\,\phi(v_{2})\,\big|
    =e^{\,\frac{1}{2}\,(u(v_{1})\,+\,u(v_{2}))}\;
    \big|\,v_{1}\,-\,v_{2}\,\big|,
  \end{equation}
  where $|\cdot|$ denotes the euclidean norm on $\R^{n}$. In other
  words, each edge length is scaled by the geometric mean of the scale
  factors $e^{u}$ attached to its vertices.

  We say that $K$ and $K'$ are \emph{discretely conformally
    equivalent} if they are discretely conformally equivalent with
  respect to some combinatorial isomorphism.
\end{definition}

This notion of discrete conformal equivalence appeared first in the
four dimensional Lorentz-geometric context of the Regge
calculus~\cite{Rocek_Williams-Quantization-1984}. In the
two-dimensional setting of surfaces, it has lead to a rich theory
which is intimately connected with hyperbolic
geometry~\cite{Bobenko_Pinkall_Springborn, Buecking-Cinfty, GuLuoWu19,
  Luo-Uniformization2, Luo-Uniformization1, Luo-Yamabe, Prosanov,
  Springborn-Uniformization} and useful in diverse applications,
see, e.g., \cite{Koehl_Hass, A01-Crane-Auxetic,
  Springborn_Schroeder_Pinkall}.

To fix ideas and introduce some notation, let us collect a few basic
facts about Möbius transformations, beginning with the definition: A
\emph{Möbius transformation} of $\widehat{\R^{n}}$ is a composition of
inversions in hyperspheres and reflections in hyperplanes, where
\begin{equation*}
  \widehat{\R^{n}}=\R^{n}\cup\{\infty\}
\end{equation*}
is the one-point compactification of $\R^{n}$. 

A Möbius transformation preserves or reverses orientation, depending
on whether it is a composition of an even or an odd number of
inversions and reflections. The Möbius transformations
of~$\widehat{\R^{n}}$ form a Lie group $\Moeb(n)$ of dimension
$\frac{1}{2}(n+1)(n+2)$ which is isomorphic to the projectivized
orthogonal group $\PO(n+1,1)$. Indeed, we may identify
$\widehat{\R^{n}}$ with the $n$-dimensional unit sphere
$S^{n}\subset\R^{n+1}$ via stereographic projection and
consider~$\R^{n+1}$ as the real projective space $\RP^{n+1}$ minus a
projective hyperplane ``at infinity''. This identifies the Möbius
group $\Moeb(n)$ with the group $\PO(n+1,1)$ of projective
transformations of $\RP^{n+1}$ that map the sphere $S^{n}$ to itself.

The group $\Sim(n)$ of similarity transformations of $\R^{n}$, i.e.,
of transformations of the form
\begin{equation*}
  x\;\longmapsto\; \lambda\, A\,x\,+\,b,
\end{equation*}
where
\begin{equation*}
  \lambda\in\R_{>0},\quad
  A\in O(n),\quad
  b\in\R^{n},
\end{equation*}
is the subgroup of Möbius transformations that fix $\infty\in\widehat{\R^{n}}$\,:
\begin{equation*}
  \Sim(n)=\{T\in\Moeb(n)\;|\;T(\infty)=\infty\}.
\end{equation*}
Conversely, the Möbius group $\Moeb(n)$ is generated by the similarity
group $\Sim(n)$ together with one sphere inversion.

In Möbius geometry, a hypersphere in $\widehat{\R^{n}}$ is either a
euclidean hypersphere in $\R^{n}$ or the union of a hyperplane in
$\R^{n}$ with $\{\infty\}$. Möbius transformations map hyperspheres to
hyperspheres. A Möbius transformation that is not a similarity
transformation does not map simplices in $\R^{n}$ to simplices, except
for zero-dimensional simplices, i.e., vertices. Two simplicial
complexes are considered Möbius equivalent if their vertices are
related by a Möbius transformation. More precisely:

\begin{definition}[Möbius equivalence]
  Simplicial complexes~$K$ and~$K'$ in~$\R^{n}$ are called
  \emph{Möbius equivalent with respect to a combinatorial isomorphism}
  $\phi:V\longrightarrow V'$ if there is a Möbius transformation
  $T\in\Moeb(n)$ such that 
  \begin{equation*}
    \phi(v)=T(v)\quad\text{for all vertices}\quad v\in V.
  \end{equation*}
  Simplicial complexes $K$ and $K'$ are called \emph{Möbius
    equivalent} if they are Möbius equivalent with respect to some
  combinatorial isomorphism.
\end{definition}

The following Definitions~\ref{def:discrete_domain}
and~\ref{def:Delaunay} specify the extra assumptions under which we
will show that discretely conformally equivalent simplicial complexes
are Möbius equivalent. All of the conditions (i)--(iii) of
Definition~\ref{def:discrete_domain} and the local Delaunay condition
of Definition~\ref{def:Delaunay} are necessary. It is easy (given the
proof of Theorem~\ref{thm:dliouville} presented in
Section~\ref{sec:proof}) to construct examples showing that the
implication may fail if any one of these conditions is not satisfied.

\begin{definition}[discrete domain]
  \label{def:discrete_domain}
  A locally finite simplicial complex $K$ in $\R^{n}$ is called a
  \emph{discrete domain} if it satisfies the following conditions:
  \begin{compactenum}[(i)]
  \item $K$ is $n$-dimensional and pure, i.e., $K$ contains only
    $n$-dimensional simplices and their faces.
  \item Every $n$-dimensional simplex in $K$ has at least one
    \emph{interior vertex}, i.e., a vertex contained in the interior
    of the carrier $|K|$.
  \item The subgraph of the $1$-skeleton of $K$ induced by the
    interior vertices is connected. 
  \end{compactenum}
\end{definition}

\begin{definition}[local Delaunay condition]
  \label{def:Delaunay}
  A discrete domain $K$ in $\R^{n}$ is called \emph{locally Delaunay}
  if, for every $n$-simplex $\sigma\in K$, the open ball bounded by
  the circumsphere of $\sigma$ contains no vertices of $n$-simplices
  sharing a common $(n-1)$-face with $\sigma$.
\end{definition}

\begin{remark}
  Let $\sigma$ and $\sigma'$ be two $n$-simplices in $\R^{n}$ that
  share a common $(n-1)$-face, say
  \begin{equation*}
    \sigma=[v_{0},\ldots,v_{n}],\qquad
    \sigma'=[v_{1},\ldots,v_{n+1}].
  \end{equation*}
  Then the following statements are equivalent:
  \begin{compactitem}
  \item $v_{0}$ is contained in the open ball bounded by the
    circumsphere of $\sigma'$.
  \item $v_{n+1}$ is contained in the open ball bounded by the
    circumsphere of $\sigma$.
  \end{compactitem}
  Thus, the local Delaunay condition imposes one condition for
  every $(n-1)$-simplex in $K$ that is incident with two $n$-simplices.
\end{remark}

To state the obvious, Theorem~\ref{thm:dliouville} refers to discrete
conformal equivalence and Möbius equivalence \emph{with respect to the
  same combinatorial isomorphism}, i.e., we will prove the following statement:

\addtocounter{theorem}{-1}
\begin{theorem}[pedantic version]
  Suppose $n\geq 3$, $K$ and $K'$ are locally Delaunay discrete
  domains in $\R^{n}$, and $\phi$ is a combinatorial isomorphism between
  $K$ and $K'$. Then $K$ and $K'$ are discretely conformally
  equivalent with respect to $\phi$ if and only if they are Möbius
  equivalent with respect to $\phi$.
\end{theorem}

\section{Proof of Theorem~\ref{thm:dliouville}}
\label{sec:proof}

\subsection{The easy implication: ``Möbius equivalence implies discrete conformal equivalence''}
\label{sec:easy}

This implication holds for arbitrary simplicial complexes in~$\R^{n}$
and for arbitrary dimension~$n$. So let $K$ and $K'$ be simplicial complexes in
$\R^{n}$ and assume they are Möbius equivalent. We will show that they
are discretely conformally equivalent.

If the Möbius transformation $T\in\Moeb(n)$ relating $K$ and
$K'$ is a similarity transformation of~$\R^{n}$, then $K$ and $K'$ are
obviously discretely conformally equivalent, because the
relation~\eqref{eq:scale} holds with a constant scale factor $e^{u}$.

If $T$ is the inversion in the unit sphere,
\begin{equation*}
  T(x)\;=\;\frac{1}{|x|^{2}}\;x,
\end{equation*}
then the identity
\begin{equation}
  \label{eq:inversion_identity}
  \big|\,
    T(x)-T(y)
  \,\big|
  \;=\;
  \frac{1}{|x|\;|y|}\;|x-y|\;,
\end{equation}
implies that $K$ and $K'$ are discretely conformally
equivalent. Indeed, in this case the relation~\eqref{eq:scale} holds
with $e^{u(v)}=|v|^{-2}$.

Since the similarity transformations and the inversion in the unit
sphere generate the Möbius group, the implication holds for all
$T\in\Moeb(n)$.

\subsection{The equivalence of simplices}
\label{sec:simplices}

In this section, we consider conformal equivalence and Möbius
equivalence for pairs of simplices. In the next section, we will use
the results to prove the harder implication, ``discrete conformal
equivalence implies Möbius equivalence.''

\begin{lemma}
  \label{lem:one_simplex}
  For $n$-simplices $[v_{0},\ldots,v_{n}]$ and
  $[v_{0}',\ldots,v_{n}']$ in $\R^{n}$, the following statements are
  equivalent:
  \begin{compactenum}[(i)]
  \item There are real numbers $u_{0},\ldots,u_{n}$ such that
    \begin{equation}
      \label{eq:simplex_scaling}
      \big|\,v_{i}'\,-\,v_{j}'\,\big|
      \;=\;
      e^{\frac{1}{2}\,(u_{i}\,+\,u_{j})}\,
      \big|\,v_{i}\,-\,v_{j}\,\big|
    \end{equation}
    for all two-element subsets $\{i,j\}$ of $\{0,\ldots,n\}$.
  \item There is a Möbius transformation\/ $T$\/ of\/
    $\widehat{\R^{n}}$ such that
    \begin{equation}
      \label{eq:T_mapping}
      v_{i}' = T(v_{i})
      \quad
      \text{for all}
      \quad
      i\in\{0,\ldots,n\}.
    \end{equation}
  \end{compactenum}
\end{lemma}

Section~\ref{sec:easy} proves the implication
``$\text{(ii)}\Rightarrow\text{(i)}$'', so it remains to show the converse
statement, ``$\text{(i)}\Rightarrow\text{(ii)}$''. This is based on the
following observations, which will also be useful by themselves:

\begin{lemma}
  \label{lem:scale}
  Assume condition (i) of Lemma~\ref{lem:one_simplex} holds. Let
  $S_{v_{0}}$ and $S_{v_{0'}}$ be the inversions in the spheres with
  radius~$1$ centered at $v_{0}$ and $v_{0}'$, respectively, and let
  \begin{equation*}
    w_{i}\;=\;S_{v_{0}}(v_{i}), \qquad w_{i}'\;=\;S_{v_{0'}}(v_{i}').
  \end{equation*}
  Then
  \begin{equation*}
    w_{0}\;=\;w_{0}'\;=\;\infty,
  \end{equation*}
  and the $(n-1)$-simplices $[w_{1},\ldots,w_{n}]$ and
  $[w_{1}',\ldots,w_{n}']$ are similar with scale factor $e^{-u_{0}}$. 
\end{lemma}

To show the implication ``$\text{(i)}\Rightarrow\text{(ii)}$'' of
Lemma~\ref{lem:one_simplex} using Lemma~\ref{lem:scale}, let~$F$ be a
similarity transformation of $\R^{n}$ mapping $w_{i}$ to $w_{i}'$ for
$i\in\{1,\ldots,n\}$, then
\begin{equation*}
  T\;=\;S_{v_{0}'}\circ F\circ S_{v_{0}}
\end{equation*}
is a Möbius transformation satisfying~\eqref{eq:T_mapping}. (Note that
there are two such similarity transformations, one of which preserves
orientation while the other reverses orientation.) This completes the
proof of Lemma~\ref{lem:one_simplex}, assuming Lemma~\ref{lem:scale}.

To prove Lemma~\ref{lem:scale}, note that
\begin{equation*}
  w_{i}\;=\;v_{0}\,+\,\frac{1}{|\,v_{i}\,-\,v_{0}\,|^{2}}\,(v_{i}\,-v\,_{0}).
\end{equation*}
for $i\in\{1,\ldots,n\}$, and a similar equation holds for $w_{i}'$.
Using the identity~\eqref{eq:inversion_identity}, one obtains
\begin{equation}
  \label{eq:wi_wj}
  |\,w_{i}\,-\,w_{j}\,|\;=\;
  \frac{1}{|\,v_{i}\,-\,v_{0}\,|\,|\,v_{j}\,-\,v_{0}\,|}
  \,|\,v_{i}\,-\,v_{j}\,|
\end{equation}
for $i,j\in\{1,\ldots,n\}$, and a similar equation for
$|\,w_{i}'\,-\,w_{j}'\,|$. Now~\eqref{eq:simplex_scaling} implies
\begin{equation*}
  |\,w_{i}'\,-\,w_{j}'\,|\;=\;e^{-u_{0}}\,|\,w_{i}\,-\,w_{j}\,|\,,
\end{equation*}
and hence the simplices are similar with scale factor $e^{-u_{0}}$.
This completes the proof of Lemma~\ref{lem:scale}.

We will also use the following fact, the proof of which we leave to
the reader:

\begin{lemma}
  \label{lem:orientation}
  If there exists any Möbius transformation $T$
  satisfying~\eqref{eq:T_mapping}, then there exist exactly two of
  them, say, $T_{1}$ and $T_{2}$, of which one preserves orientation
  while the other reverses orientation, and which are related by
  \begin{equation*}
    T_{2}= T_{1} \circ C = C' \circ T_{1},
  \end{equation*}
  where $C$ and $C'$ are the inversions in the circumspheres of the
  simplices~$[v_{0},\ldots,v_{n}]$ and~$[v_{0}',\ldots,v_{n}']$,
  respectively.
\end{lemma}

\subsection{The harder implication: ``Discrete conformal equivalence
  implies Möbius equivalence''}

\label{sec:harder}

Let $K$ and $K'$ be two locally Delaunay discrete domains in
$\R^{n}$, where $n\geq 3$, and assume $K$ and $K'$ are discretely
conformally equivalent with respect to the combinatorial isomorphism
$\phi$.

Note that $\phi$ may be orientation preserving or orientation
reversing. We may assume without loss of generality that $\phi$ is
orientation preserving. (If $\phi$ is orientation reversing, consider
orientation preserving isomorphism between $K$ and a mirror image of
$K'$.)

Lemmas~\ref{lem:one_simplex} and~\ref{lem:orientation} say that for
each $n$-simplex $\sigma\in K$, there is a unique orientation
preserving Möbius transformation $T_{\sigma}$ such that
$T_{\sigma}(v)=\phi(v)$ for every vertex $v\in\sigma$. Note that the
assumptions about orientation ensure that $T_{\sigma}$ maps the inside
of the circumsphere of $\sigma\in K$ to the inside of the circumsphere
of $\phi(\sigma)\in K'$.

It remains to show that these Möbius transformations $T_{\sigma}$ are
in fact all equal. To this end, it is enough to show the equality for
$n$-simplices in the star of an interior vertex, i.e., to
show the following lemma:

\begin{lemma}
  \label{lem:star}
  If $v$ is an interior vertex of $K$ and if $\sigma$ and
  $\tilde\sigma$ are two $n$-simplices contained in $\st(v)$, then
  $T_{\sigma}=T_{\tilde{\sigma}}$.
\end{lemma}

Indeed, suppose Lemma~\ref{lem:star} holds and $\sigma$ and
$\tilde{\sigma}$ are any $n$-simplices of~$K$. By assumption, both
contain interior vertices of $K$, say $v$ and
$\tilde{v}$. Furthermore, by assumption, there is a path from~$v$
to~$\tilde{v}$ in the $1$-skeleton of $K$ traversing only interior
vertices. By induction on the length of the path, Lemma~\ref{lem:star}
implies that $T_{\sigma}=T_{\tilde{\sigma}}$.

The proof of Lemma~\ref{lem:star} relies on Cauchy's rigidity theorem
for convex polyhedra, applied to the link of $v$ after an inversion
centered at $v$. Convexity follows from the following general
observation:

\begin{lemma}\label{lem:Delaunay-convex}
  Let $K_{0}$ be a simplicial complex in $\R^n$, $n\geq 1$, with one
  interior vertex $v$ and $K_{0}=\st(v)$. Let $S$ be an inversion in
  some sphere centered at $v$, and apply it to the vertices of
  $\partial K_{0}$ to obtain a Möbius equivalent $(n-1)$-dimensional
  complex $P$. Then the following statements are equivalent:
  \begin{compactenum}[(i)]
  \item $K_{0}$ satisfies the local Delaunay condition.
  \item $P$ is convex.
  \end{compactenum}
\end{lemma}

\begin{proof}[Proof (Lemma~\ref{lem:Delaunay-convex})]
  Each $n$-simplex of $K_{0}$ corresponds to an
  $(n-1)$-di\-men\-sio\-nal face of~$P$, and $S$ maps the closed ball
  bounded by the circumsphere of an $n$-simplex of~$K_{0}$ to a closed
  halfspace in $\R^{n}$ bounded by the hyperplane of the respective
  face of~$P$. Thus, the local Delaunay condition for $K_{0}$ is
  equivalent to a local convexity condition for $P$ involving adjacent
  faces, which is equivalent to the global convexity of~$P$.
\end{proof}

Now to show Lemma~\ref{lem:star}, let
$v'=\phi(v)$, and let $Q$ and $Q'$ be the links of~$v$ and~$v'$,
respectively, i.e.,
\begin{equation*}
  Q=\partial\st(v),\qquad Q'=\partial\st(v').
\end{equation*}
Let $S_{v}$ and $S_{v'}$ be the inversions in the spheres with radius
$1$ centered at~$v$ and~$v'$, respectively. Apply the inversions
$S_{v}$ and $S_{v'}$ to the vertices of $Q$ and $Q'$, respectively, to
obtain for each a Möbius equivalent polyhedron, $P$ and $P'$. By
Lemma~\ref{lem:Delaunay-convex}, the local Delaunay condition for $K$
and $K'$ implies that $P$ and $P'$ are convex polyhedra in
$\R^{n}$. As in Section~\ref{sec:simplices}, one sees that the facets
of $P$ and $P'$ are similar. By Cauchy's rigidity theorem for convex
polyhedra and its higher dimensional generalization~\cite{Pak}, there
is a similarity transformation $F$ of $\R^{n}$ that maps $P$ to
$P'$. By the orientation assumption, $F$ is orientation
preserving. Hence
\begin{equation*}
  T = S_{v'}\circ F\circ S_{v}
\end{equation*}
is an orientation preserving Möbius transformation that maps $\st(v)$
to $\st(v')$. Therefore, $T=T_{\sigma}$ for all
$\sigma\in\st(v)$. This proves Lemma~\ref{lem:star} and hence the
implication ``discrete conformal equivalence implies Möbius
equivalence'', and this completes the proof of
  Theorem~\ref{thm:dliouville}.

\section{Discrete conformal flatness and induced hyperbolic metric:
  Concluding remarks, outlook and open questions}
\label{sec:outlook}

The definition of discrete conformal equivalence
(Definition~\ref{def:dce}) extends in an obvious way from simplicial
complexes in $\R^{n}$ to triangulated piecewise euclidean manifolds,
possibly with boundary, i.e., to manifolds that consist of euclidean
simplices glued together along their facets. We propose the following
notion of discrete conformal flatness:

\begin{definition}
  A triangulated piecewise euclidean manifold is \emph{discretely
    conformally flat} if the vertex star of every interior vertex is
  discretely conformally equivalent to a locally Delaunay discrete
  domain.
\end{definition}

Applying the same idea as in the proof of the discrete Liouville
theorem, and in particular equation~\eqref{eq:wi_wj}, one obtains the
following result:

\begin{theorem}
  A $n$-dimensional triangulated piecewise euclidean manifold is
  discretely conformally flat if and only if every interior vertex
  $v_{0}$ satisfies the following condition:

  There exists a convex polyhedron in $\R^{n}$ that is combinatorially
  equivalent to the link of $v_{0}$ and whose edge lengths are
  \begin{equation}
    \tilde{\ell}_{ij}=\frac{\ell_{ij}}{\ell_{0i}\ell_{0j}}\;,
  \end{equation}
  Here, $\ell_{ij}$ denotes the length of the edge between two
  adjacent vertices $v_{i}$, $v_{j}$ of the triangulated piecewise
  euclidean manifold, and for vertices $v_{i}$, $v_{j}$ in the
  link of $v_{0}$, $\tilde{\ell}_{ij}$ denotes length of the
  corresponding edge of the convex polyhedron.
\end{theorem}

Note that for $n=2$, the condition on the vertex link is equivalent to
the \emph{polyhedral inequalities} for the $\tilde{\ell}_{ij}$, i.e.,
each $\tilde{\ell}_{ij}$ is larger than the sum of the others.

Note also that a connection between discrete conformal equivalence and
hyperbolic geometry, which plays an important role in the theory in
dimension two~\cite{Bobenko_Pinkall_Springborn,
  Springborn-Uniformization}, extends to higher dimensions: If you
interpret the circumsphere of a euclidean $n$-simplex as the boundary
of $n$-dimensional hyperbolic space in the Beltrami-Klein model, this
induces a hyperbolic metric on the simplex minus its vertices, turning
the simplex into an ideal hyperbolic simplex. If you perform this
construction on all simplices of a triangulated piecewise euclidean
manifold, this induces a hyperbolic metric on the manifold with cusps
at the vertices and cone-like singularities in the faces of
codimension two. Just as in the two-dimensional setting, one can prove
the following theorem:

\begin{theorem}
  Two triangulated piecewise euclidean manifolds are discretely
  conformally equivalent if and only if they are isometric with
  respect to the induced hyperbolic metrics. 
\end{theorem}

As in the two-dimensional setting~\cite{Bobenko_Pinkall_Springborn,
  Springborn-Uniformization}, this observation suggests extending the
definition of discrete conformal equivalence to triangulations that
are not combinatorially equivalent:

\begin{definition}[discrete conformal equivalence, extended] Two
  triangulated piecewise euclidean manifolds (which need not be
  combinatorially equivalent) are \emph{discretely
    conformally equivalent} if they are isometric with respect to the
  induced hyperbolic metrics.
\end{definition}

In the $2$-dimensional setting, the known uniformization
results~\cite{Luo-Uniformization1,Springborn-Uniformization} show that
any triangulated surface is discretely conformally flat, provided the
notion of discrete conformal flatness is based on the extended notion
of discrete conformal equivalence.

In dimensions $3$ and higher, the situation more complicated. The
induced hyperbolic metric will in general have cone-like singularities
along faces of codimension $2$, even if the piecewise euclidean
manifold is flat. The total dihedral angles at the faces of
codimension $2$ are discrete conformal invariants. Thus, any
codimension-$2$-face whose hyperbolic cone angle is not equal to
$2\pi$ occurs in any discretely conformally equivalent triangulated
piecewise euclidean manifold.

Maybe the most intriguing question opened by this line of research is
how far the analogy between smooth and discrete conformal flatness
extends. Consider for example the case of closed $3$-dimensional
manifolds $M$. The Chern--Simons functional $CS_{M}(g)$ is an
$\R/\Z$-valued conformally invariant function on the space of
Riemannian metrics $g$ on $M$, and the critical points of $CS_{M}$ are
precisely the conformally flat metrics on $M$, see, e.g.,~\cite{Moroianu}. Is
there an analogous invariant on the set of discrete conformal classes
whose critical points are exactly the discretely conformally flat
classes?

\section*{Acknowledgments}

This research was supported by the DFG SFB/TR 109 ``Discretization in Geometry
and Dynamics''.

\begingroup
\small
\bibliographystyle{abbrv}
\bibliography{liouville}

\bigskip\noindent%
\textit{Technische Universität Berlin,
  Institut für Mathematik,
  Strasse des 17. Juni 135,
  10623 Berlin, Germany}

\medskip\noindent%
\href{mailto:pinkall@math.tu-berlin.de}{\nolinkurl{pinkall@math.tu-berlin.de}}\\
\href{mailto:boris.springborn@tu-berlin.de}{\nolinkurl{boris.springborn@tu-berlin.de}}

\medskip\noindent%
\url{https://page.math.tu-berlin.de/~pinkall}\\
\url{https://page.math.tu-berlin.de/~springb}

\endgroup

\end{document}